\documentclass{amsart}
\usepackage{graphicx}
\theoremstyle{definition}

\theoremstyle{remark}

\numberwithin{equation}{section}
\begin{document}
\title[]{On the Ergodic Properties of Certain Additive Cellular Automata over $Z_{m}$}%
\author{Hasan Akin}%
\address{}%
\email{}%

\thanks{}%
\subjclass{Primary 28D20; Secondary 37A15 }%
\keywords{Cellular Automata, ergodicity, weak-mixing, strong mixing}%

%\date{}%
%\dedicatory{}%
%\commby{}%
% ----------------------------------------------------------------
\begin{abstract}
In this paper, we investigate some ergodic properties of
$Z^{2}$-actions $T_{p,n}$ generated by an additive cellular
automata and shift acting on the space of all doubly -infinitive
sequences taking values in $Z_{m}$.
\end{abstract}
\maketitle \maketitle
% ----------------------------------------------------------------
\section{Introduction}
Mathematical study of cellular automata was initiated by Hedlund
late 1960s. Hedlund determined the properties of endomorphisms and
automorphisms of the shift dynamical system[2]. Sato studied
linear cellular automata with-dimensional cell space as well as
higher-dimensional cell space[3]. The properties of endomorphisms
of subshifts of finite type were studied by Coven et al. [1].
Sinai gave a formula for directional entropy[5]. Ergodic
properties of cellular automata have been investigated from
various aspects by Shereshevsky and proved that if the automata
map is bipermutative then associated CA- action is
strongly-mixing[4].

\par In this paper, we shall restrict our attention to additive
cellular automata over $Z_{m}$. The organization of the paper is
as follows: In section 2 we establish the basic formulation of
problem necessary to state our main theorem. In section 3 we prove
our main theorem and some results. Let us provide some notation
and background.

\section{Formulation of the problem }
Let \ $Z_{m}=\left\{ 0,1,...,m-1\right\} $ be a finite alphabet
and $\Omega =Z_{m}^{Z}$ be the space of double-infinite sequences
$x=\left( x_{n}\right) _{n=-\infty }^{\infty }$, $x_{n}\in Z_{m}$,
$\sigma$  is the shift in $\Omega$, i.e. $\sigma
x=x^{^{\prime}}=\left\{ x_{n}^{^{\prime }}\right\} $,
$x_{n}^{^{\prime }}=x_{n+1}$, $x_{n}\in Z_{m}$. A continuous map
$f_{\infty }:\Omega \rightarrow \Omega $  commuting with the shift
(i.e. such that $f_{\infty }\circ \sigma =\sigma \circ f_{\infty
}$) is called a cellular automaton. It is well known (see([2],
Theorem 3.4)) that $f_{\infty }:\Omega \rightarrow \Omega $ is a
cellular automaton if and only if there exist $l,r\in Z$ with
$l\leq r$ and a mapping $f:Z_{m}^{r-l+1}\rightarrow Z_{m}$ such
that
\begin{eqnarray*}
f_{\infty }(x)=(y_{n})_{n=-\infty }^{\infty },
y_{n}=f(x_{n+l},...,x_{n+r})
\end{eqnarray*}
for all $x\in \Omega $. $n\in Z$. It is called the mapping $f$ the
rule of  $f_{\infty }$ and the interval $\left[ l,r\right] $ the
range of $f_{\infty }$. In [5], it was assumed that $\sigma$ and
$f_{\infty }$ generate an action of the group $Z^{2}$ on $\Omega
$: for $(m, n)\in Z^{2}$ the corresponding transformation is
$T_{p,n}=\sigma ^{p}f_{\infty }^{n}$. Firstly, we consider
additive cellular automata $f_{\infty }$determined by an
automation rule
\begin{eqnarray*}f(x_{n-k},...,x_{n+k})=
(\overset{k}{\underset{i=-k}{\sum }}\lambda
_{i}x_{n+i})({mod}m)(\lambda _{i}\in Z_{m}).
\end{eqnarray*}
A cellular automaton (CA) defined on $\Omega$ is a map $F:\Omega
\rightarrow \Omega$ such that for $x \in \Omega$ and $i \in Z$,
$(Fx)_{i}=f(x_{i-r},...,x_{i+r}$ where $r \in N$ is radius and
$f:\mathbb{Z}_{m}^{2r+1}\rightarrow \mathbb{Z}_{m}$ is a given
local rule. Generally, we take as $(\lambda _{i}=1)$. Let us
consider a block $A=_{a-k}\left[ i_{a-k},...,i_{a+k}\right]
_{a+k}$. The first preimage of the block $A$ under $f_{\infty }$
is
\\
$\left\{ y\in \Omega :y_{a-2k}=j_{a-2k},...,y_{a+2k}=j_{a+2k},
j_{a-2k},...,j_{a+2k}\in Z_{m}\right\} $ where

$y_{a-2k}+...+y_{a}=i_{a-k}({mod}m)$,

  ....

  ....

  ....

$y_{a-k}+...+y_{a+k}=i_{a}({mod}m)$,

  ....

  ....

  ....

$y_{a}+...+y_{a+2k}=i_{a+k}({mod}m)$.
\\

\par It is easy to see from this system of equations that $(f_{\infty })^{-1}(A)$
consists of $m^{2k}$ following blocks $(j_{a-2k},...,$
$j_{a+2k})$. Now we calculate the measure

\begin{eqnarray*}
\mu((f_{\infty})^{-1}(A))& = & m^{2k}\mu \{y\in \Omega
:y_{a-2k}=j_{a-2k},...,y_{a+2k}=j_{a+2k}, j_{a-2k},
j_{a+2k}\in Z_{m}\}\\
& = & m^{2k}m^{-(4k+1)}=m^{-(2k+1)}.
\end{eqnarray*}

\par {\bf Example.} Let $A=\left\{ 0,1\right\} $ and
$f(x_{-2},x_{-1},x_{0},x_{1},x_{2})=\left(
\underset{i=-2}{\overset{2} {\sum }}x_{i}\right) $(mod $2$). Then
\\

$\left( f_{\infty }\sigma \right) ^{-1}\left( _{-2}
\left[10101\right] _{2}\right) =_{\text{ }-3} \left[
111110000\right] _{5}\cup $ $_{-3}\left[ 100000111\right] _{5}\cup
$ $_{-3}\left[ 010001011\right] _{5}$
\\

\qquad\qquad\qquad\qquad\qquad $\cup _{\text{ }-3}\left[
001001101\right] _{5}\cup _{\text{ }-3}\left[ 000101110\right]
_{5}\cup _{\text{ }-3}\left[ 000011111\right] _{5}$
\\

\qquad\qquad\qquad\qquad\qquad $\cup _{\text{ }-3}\left[
111000001\right] _{5}\cup _{\text{ }-3}\left[ 011101000\right]
_{5}\cup _{\text{ }-3}\left[ 001111100\right] _{5}$
\\

\qquad\qquad\qquad\qquad\qquad $\cup _{\text{ }-3}\left[
110100010\right] _{5}\cup _{\text{ }-3}\left[ 110010011\right]
_{5}\cup _{\text{ }-3}\left[ 101100100\right] _{5}$
\\

\qquad\qquad\qquad\qquad\qquad $\cup _{\text{ }-3}\left[
100110110\right] _{5}\cup _{\text{ }-3}\left[ 101010101\right]
_{5}\cup _{\text{ }-3}\left[ 010111010\right] _{5}$
\\

\qquad\qquad\qquad\qquad\qquad$\cup _{\text{ }-3}\left[
011011001\right] _{5}$. Thus we have
\begin{eqnarray*}
\mu(\left( f_{\infty }\sigma \right)
^{-1}(_{-2}[10101]_{2}))=16\mu
(_{-3}[j_{-4},...,j_{4}]_{5})=2^{4}2^{-9}=2^{-5}.
\end{eqnarray*}
If we continue this operation, by the same way, we can determine
the measure of (n-1)st preimage of the block
$A=_{a-k}[i_{a-k},...,i_{a+k}]_{a+k}$ under $f_{\infty }$.
\par Evidently this (n-1)st preimage consist of such $(z_{n})_{n=-\infty }^{\infty }$, for
which we have following system of equations:

$z_{a-nk}+...+z_{a-(n-1)k}+...+z_{a-(n-2)k}=h_{a-(n-1)k}({mod}m)$,

  ....

  ....

  ....

$z_{a-k}+...+z_{a}+...+z_{a+k}=h_{a}({mod}m)$,

  ....

  ....

  ....

$z_{a+(n-2)k}+...+z_{a+(n-1)k}+...+z_{a+nk}=h_{a+(n-1)k}({mod}m)$,
\\
where $h_{a-(n-1)k},...,h_{a},...,h_{a+(n-1)k}\in Z_{m}$. So we
can calculate the measure
\begin{eqnarray*}
\mu(f_{\infty}^{-(n-1)}(A))=m^{2(n-1)k}m^{-(2nk+1)}=m^{-(2k+1)}.
\end{eqnarray*}
\section{Results}
Here we shall use the terminology of Sinai [5]. Let us consider as
$Z^{2}-action$ $T_{p,n}=\sigma ^{p}f_{\infty }^{n}$.
\\

\par {\bf Proposition:} Let $T_{p,n}=\sigma ^{p}f_{\infty
}^{n}$ be $Z^{2}-action$ as above and if $\mu$ is stationary
Bernoulli measure on $\Omega$, that is, $\mu (i)=\frac{1}{m},$
$\forall \ i=0,1,...,m-1$, then both $f_{\infty }$ and $T_{p,n}$
are Bernoulli measure preserving transformations.
\\

\par {\bf Lemma:} The surjective CA-map $f_{\infty }$
generated by the rule
\begin{center}
$f(x_{n+l},...,x_{n+r})=(\overset{r}{\underset{i=l}{\sum
}}x_{n+i})(modm)$
\end{center}
is nonergodic with respect to the measure $\mu$, because the
equality
$$
\mu(_{b}[e_{0},...,e_{s}]_{b+s}\cap
f_{\infty}^{-n}(_{a}[d_{0},...,d_{k}]_{a+k}))
 = \mu(_{b}[e_{0},...,e_{s}]_{b+s})\mu(_{a}[d_{0},...,d_{k}]_{a+k})
$$
can't be obtained sometimes. But we show that $Z^{2}-action$
$T_{p,n}=\sigma^{p}f_{\infty }^{n}$ defined by $(p,n)\mapsto
T_{p,n}=\sigma^{p}f_{\infty }^{n}$ on $(\Omega, \mathcal{B}, \mu)$
is ergodic, weak-mixing and strong-mixing if $p>b+s+n\ell-a$.
\\

\par {\bf Theorem 1:} [6, Theorem 1.17] Let (X,$\mathcal{B}$,
$\mu$) be a measure space and let $\mathcal{A}$ be a semi-algebra
that generates $\mathcal{B}$. Let T:X $\rightarrow$ X be a
measure-preserving transformation. Then
\par(i) T is ergodic iff $\forall A,B\in$$\mathcal{A}$
\begin{eqnarray*}\underset{n\rightarrow \infty
}{\lim}\frac{1}{n}\overset{n-1}{\underset{i=0}{\sum }}\mu
(T^{-i}A\cap B)=\mu (A)\mu (B),
\end{eqnarray*}
\par(ii) T is weak-mixing iff $\forall$ A,B$\in$$\mathcal{A}$
\begin{eqnarray*}\underset{n\rightarrow \infty }{\lim }\overset{n-1}
{\underset{i=0}{\sum }}\left| \mu (T^{-i}A\cap B)-\mu (A)\mu
(B)\right| =0
\end{eqnarray*} and
\par(iii) T is strongly-mixing iff $\forall A,B\in $ $\mathcal{A}$
\begin{eqnarray*}
\underset{n\rightarrow \infty }{\lim }\mu (T^{-n}A\cap B)=
\mu(A)\mu(B).
\end{eqnarray*}
Now we can give main theorem.
\par {\bf Theorem 2:} Let $Z_{m}=\left\{0,1,...,m-1\right\} $
be a finite alphabet and $\Omega =Z_{m}^{Z}$ be the space of
double-infinite sequences
 $x=\left( x_{n}\right) _{n=-\infty }^{\infty }$, $x_{n}\in
Z_{m}$. If additive cellular automata $f_{\infty }$ is given by
the formula:
\begin{eqnarray*}f_{\infty }(x)=(y_{n})_{n=-\infty
}^{\infty },
y_{n}=f(x_{n+\ell},...,x_{n+r})=(\overset{r}{\underset{i=\ell}{\sum
}}x_{n+i})({mod}m)
\end{eqnarray*}
 for all $x\in \Omega $. $(p,n)\in
Z^{+} \times Z^{+}$, then $Z^{2}-action$ $T_{p,n}=\sigma
^{p}f_{\infty }^{n}$ is ergodic, strongly-mixing and weak-mixing.

\begin{proof} To prove that $T_{p,n}$ is ergodic it is sufficient
to verify (Theorem 1,ii)for any two cylinder sets
$A=_{a}[d_{0},...,d_{k}]_{a+k}$ and
$B=_{b}[e_{0},...,e_{s}]_{b+s}$, we have
\begin{eqnarray*} \underset{p,n\rightarrow \infty }{\lim
}\frac{1}{pn}\underset{(i,j)\in D}{\sum }\mu
(_{b}[e_{0},...,e_{s}]_{b+s}\cap
T_{({-i,-j})}(_{a}[d_{0},...,d_{k}]_{a+k}))=
\\
\mu(_{b}[e_{0},...,e_{s}]_{b+s})\mu (_{a}[d_{0},...,d_{k}]_{a+k}),
\end{eqnarray*}
where $D=[0,p-1]\times[0,n-1] \cap Z^{2}$. For i$>$b+s+j$\ell$-a
we have
\begin{eqnarray*}
\mu (_{b}[e_{0},...,e_{s}]_{b+s}\cap
T_{(-i,-j)}(_{a}[d_{0},...,d_{k}]_{a+k}))
 & = &
 \mu(_{b}[e_{0},...,e_{s}]_{b+s})\mu
(_{a}[d_{0},...,d_{k}]_{a+k}).
\end{eqnarray*}
On the other hand, we show that
\begin{eqnarray*} \underset{p,n\rightarrow \infty
}{\lim}\frac{1}{pn}\underset{(i,j)\in D}{\sum
}\mu(_{b}[e_{0},...,e_{s}]_{b+s}\cap
T_{(-i,-j)}(_{a}[d_{0},...,d_{k}]_{a+k}))&&\\
\end{eqnarray*}
\begin{eqnarray*}
&  =  &
 \underset{p,n\rightarrow \infty }{\lim }\frac{1}{pn}\mu
(_{b}[e_{0},...,e_{s}]_{b+s})\underset{(i,j)\in D}{\sum }f_{\infty
}^{-j}\sigma ^{-i}(_{a}[d_{0},...,d_{k}]_{a+k}))\\
& = & \mu(_{b}[e_{0},...,e_{s}]_{b+s})\underset{p,n\rightarrow
\infty }{\lim }\frac{1}{pn}\underset{(i,j)\in D}{\sum }f_{\infty
}^{-j}(_{a+i) }[d_{0},...,d_{k}]_{a+k+i}))\\
& = &
 \mu (B)\underset{p,n\rightarrow \infty }{\lim
}\frac{1}{pn}\overset{n-1}{\underset{j=0}{\sum }}(\mu (f_{\infty
}^{-j}(_{a}[d_{0},...,d_{k}]_{a+k})+... +\mu (f_{\infty
}^{-j}(_{a+p-1}[d_{0},...,d_{k}]_{a+k+p-1}))\\
& = & \mu(B)\underset{p,n\rightarrow \infty }{\lim
}\frac{1}{pn}\overset{n-1}{\underset{i=0}{\sum }}[pm^{-(k+1)}]\\
& = & \mu (B)\mu (A).
\end{eqnarray*}
\\
So $Z^{2}-action$ $T_{p,n}=\sigma ^{p}f_{\infty }^{n}$ is ergodic.
Similarly for $i>b+s+j\ell-a$ we have
\begin{eqnarray*}
\mu(_{b}[e_{0},...,e_{s}]_{b+s}\cap
T_{(-i,-j)}(_{a}[d_{0},...,d_{k}]_{a+k})) =
\mu(_{b}[e_{0},...,e_{s}]_{b+s})\mu (_{a}[d_{0},...,d_{k}]_{a+k}).
\end{eqnarray*}
Let $A=_{a}[d_{0},...,d_{k}]_{a+k}$ and
$B=_{b}[e_{0},...,e_{s}]_{b+s}$ be any arbitrary cylinder sets.
Then we have
\begin{eqnarray*}
\underset{p,n\rightarrow \infty }{\lim }\mu [T_{(-p,-n)}(A)\cap
B]&=& \underset{p,n\rightarrow \infty }{\lim }\mu
[(f_{\infty})^{-n}(_{a+p}[d_{0},...,d_{k}]_{a+k+p}\cap B]\\
 &  = &
\mu (B)\underset{p,n\rightarrow \infty }{\lim }(\mu f_{\infty
}^{-n}(_{a+p}[d_{0},...,d_{k}]_{a+k+p}))\\
&  = & \mu (B)\mu (A).
\end{eqnarray*}
Because every strongly-mixing transformation is weak-mixing,
$T_{(p,n)}$ is weak-mixing.
\end{proof}
\par \par One can prove that the natural extension of $T_{p,n}=\sigma ^{p}f_{\infty }^{n}$
is ergodic and mixing.
\begin{center}Acknowledgement
\end{center}

\par The author is grateful to Professor Nasir Ganikhodjaev for
encouragement and many helpful discussions during the preparation
of this work.

\vspace{1cm}

\baselineskip=0.4\baselineskip

{\small
\noindent Harran University \\
Arts and Sciences Faculty \\
Department of Mathematics \\
 6300, \c Sanl\i urfa, TURKEY \\
e-mail: akinhasan@harran.edu.tr


\begin{thebibliography}{6}
    \bibitem{1} E. M. Coven and M. E. Paul, Endomorphisms of irreducible subshift of finite
    type, Math. Sys. Theory, \textbf{8} (1974), 167-175.
    \bibitem{2} G. A. Hedlund, Endomorphisms and automorphisms of the shift dynamical system,
     Math. Sys. Theory, \textbf{3} (1969), 320-375.
    \bibitem{3} T. Sato, Ergodicity of Linear Cellular Automata over $Z_{m}$,
    Inform.Processing Letters \textbf{61} (1997), 169-172.
    \bibitem{4} M. A. Shereshevsky, Ergodic properties of certain surjective cellular automata,
     Monatsh. Math. \textbf{114} (1992), 305-316.
    \bibitem{5} Ya. G. Sinai, An answer to a question by J. Milnor,
    Comment. Math. Helvetici \textbf{60} (1985), 173-178.
    \bibitem{6} P. Walters, An Introduction to Ergodic Theory,
    Springer-Verlag (1982).
\end{thebibliography}
\end{document}